\documentclass[notitlepage]{article}
\usepackage{amssymb,amsmath,amsfonts}

\newif\ifcomments
\commentsfalse

\newcommand\unmig{\frac12}

\newcommand\tp{{\theta_p}}
\newcommand\tpu{{\theta_{p+1}}}

\newcommand\oeth{\hbox{\lower.3ex\hbox{\large l}\kern-.2em\raise.1ex\hbox{o}}}

\newtheorem{theorem}{Theorem}

\newtheorem{lemma}[theorem]{Lemma}

\newcommand{\demo}[1]{\par\smallskip \noindent {\bf {#1}:}}
\newcommand{\qed}{\ \hfill\mbox{$\clubsuit$}}

\def\begeq{\begin{equation}}
\newenvironment{comment}{\small \quotation }{\endquotation\normalsize}

\ifcomments
\def\nota#1{\begin{comment}{#1}\end{comment}}
\else
\def\nota#1{\relax}
\fi

\title{Yet another inverse function theorem}

\author{
Jaume Gudayol \\
Dept. de Matem{\`a}tica Aplicada i An{\`a}lisi,\\
 Universitat de Barcelona. \\
 E-mail:{\it gudayol@mat.ub.es}
}

\date{September 29, 1992}

\begin{document}

\newpage

\begin{flushright}
\begin{minipage}{10cm}
\begin{flushright}
   \ldots \ \ \ Ic to so{\oeth}e wat \\
\hfill \hbox{{\oeth}{\ae}{}t bi{\oeth} in eorle  \ \  indryhten {\oeth}eaw,}\\
\hfill \hbox{{\oeth}{\ae}{}t he his fer$\eth$locan \ \ \  f{\ae}{}ste binde,}\\
healde his hordcofan, \ \ \  hycge swa he wille. \\

The Wanderer. 
\end{flushright}
\end{minipage}
\end{flushright}

\let\prova=\newpage

\let\newpage=\relax

\maketitle

\let\newpage=\prova

\begin{abstract}
We prove a Nash-Moser type inverse function theorem in Fr{\'e}chet
spaces for functions with approximate inverses, allowing for a loss of
derivatives proportional to $n$ in the way of Lojasiewicz and Zehnder.
\end{abstract}

\subsection{Introduction}
\nota{
This article was written (in catalan) as an end of term exercise for
the course `Non linear fuctional analysis', given by Ernest Fontich
in the 1991-1992 term. Since I have been told that it has circulated 
since then, I have decided to make it public.
}
The goal of this paper is to obtain an inverse function theorem for
functions with approximative inverses in graded Fr{\'e}chet spaces,
similar to the one that can be 
found in \cite{Zehnder}. There, a Nash-Moser inverse function theorem
is proven for functions $\phi$ such that its derivative $D\phi$ does not
have a right inverse, but only an approximation to it. In
\cite{Zehnder}, such a theorem is proven, under the assumption that 
the loss of derivatives does not depend on $n$. Thanks to Lojasiewicz
and Zehnder (see \cite{L-Z}), it is known that, when a right inverse
for $D\phi$ exists, the loss of derivatives can be of the order of
$(\lambda-1) n$, as far as $\lambda<2$. Our goal is to provide a bridge
between these two theorems. 

All along the article, $E$ and $F$ will be tame Fr\'echet spaces with 
graded norms, and $S_\theta$, for $\theta\ge1$, will be the
corresponding smoothing operators. The map $\phi$ of which we want to
find the inverse will be a tame map, with tame derivative, that is,
$|\phi(x)|_n\le C|x|_n$ and 
$|\phi'(x)v|_n\leq C\left(|x|_{n}|v|_d+|v|_{n}\right)$.
For such a function, we prove the following
\setcounter{theorem}{-1}
\begin{theorem}
Let $E$, $F$, and $\phi:E\mapsto F$ be as above. Assume that there is
a local approximate right inverse $L:(U\subset E)\times F\to E$ to
$\phi'$, where 
$U$ is some neighborhood of $0\in E$, satisfying
$$
|(\phi'(x)L(x)-I)y|_n\leq C
\left(|x|_{n}|y|_d+|y|_{n}\right) |x|_n
$$
and
$$
|(\phi'(x)L(x)-I)y|_n\leq C
\left(|x|_{n}|y|_d+|y|_{n}\right),
$$
for some $d$. Assume also that 
$$
|L(x)y|_n\leq C\left(|x|_{\lambda n+d}|y|_d+|y|_{\lambda n
+d}\right),
$$
for any $(x,y)\in U\times F$ and $n\geq0$, and given $d\geq0$ and
$1\le\lambda<2$. Assume also that for the Taylor rest 
$R(x,v)=\phi(x+v)-\phi(x)-\phi'(x)v$ there are $C$ and $l$ so that 
$$
|R(x,v)|_n\leq C\left(|x|_n|v|_l^2+|v|_l|v|_n\right).
$$
Moreover, assume that there is an $m$ so that
$$
|\sum_{l=0}^{\infty}(S_\theta(1-\phi'(x)L(x)))^l y|_n\leq
C(n)\theta^{m} (|x|_n |y|_d+ |y|_n).
$$
Then $\phi$ has a local right inverse $\psi$,
defined on a neighborhood $V=\{y\in F \,|\,|y|_{s_0}<\delta\}$ of $0$, and
satisfying $|\psi(y)|_d\leq C|y|_{s_0}$, for some $s_0$. 
\end{theorem}


\subsection{The spaces}
Let $E$ and $F$ be two Fr{\'e}chet spaces, with graded norms (for the
definitions, see Hamilton \cite{Hamilton}).
Moreover, we will assume that $E$ and $F$ are tame in the sense that
there exist in each of them a collection of smoothing linear operators 
$\{S_\theta\,|\,\theta\geq1\}$, from the space into itself, such that,
for $0\leq k\leq n$: 
\begin{eqnarray}
 |S_\theta x |_n  &\leq & C(n,k)\theta^{n-k}|x|_k \notag \\
         |(1-S_\theta)x|_k &\leq & C(n,k)\theta^{-(n-k)}|x|_n
\label{i}
\end{eqnarray}
Recall that from \ref{i} we get the interpolation inequalities:
$$
|x|_l\leq C(l,k,n)\,|x|_k^{1-\alpha}|x|_n^\alpha\qquad\mbox{ if }
l=(1-\alpha)k+\alpha n.
$$

\subsection{The function}
Let $\phi:E\to F$ be a continuous function such that $\phi(0)=0$. 
We want to know which conditions on $\phi$ assure us the existence of
a local inverse $\psi$, defined on a neighborhood $V$ of $0\in F$ and
satisfying  $\phi\circ\psi=I$ on $V$. The conditions we will impose on
$F$ will apply in a neighborhood of $U=\{|x|_l<1\}$, for some given $l>0$.
The first condition we will require is that $\phi$ be tame, that is
(using that $\phi(0)=0$) that there is a $C$ for which: 
$$
|\phi(x)|_n\leq C|x|_{n+d_1},
$$
whenever $x\in U$, for all $n\geq 0$ and a given $d_1$.
Let us notice, though, that by renumbering the seminorms, we can get
$d_1=0$, and thus the requirement is:
\begin{equation}
|\phi(x)|_n\leq C|x|_{n}.\label{1}
\end{equation}
Another requirement is that $\phi$ be Fr{\'e}chet differentiable in
$U$, and that its dedivative can be bounded by
$$
|\phi'(x)v|_n\leq C\left(|x|_{n+d_2}|v|_l+|v|_{n+d_2}\right)
$$
for some $d_2\geq0$, and all $(x,v)\in U\times E$ and $n\geq0$. Again,
by reordering the seminorms we can get $d_2=0$. Thus what we ask for
is that, for some $C$ and $d$:
\begin{equation}
|\phi'(x)v|_n\leq C\left(|x|_{n}|v|_d+|v|_{n}\right).\label{2}
\end{equation}
Moreover, we want $\phi'$ to have an approximate right inverse, that
is, we ask for the existence of a function $L:(U\subset E)\times F\to E$,
lineal with respect to $F$, satisfying:
\begin{equation}
|(\phi'(x)L(x)-I)y|_n\leq C
\left(|x|_{n}|y|_d+|y|_{n}\right)
   |x|_n\label{3}
\end{equation}
and
\begin{equation}
|(\phi'(x)L(x)-I)y|_n\leq C
\left(|x|_{n}|y|_d+|y|_{n}\right)\label{4}
\end{equation}
In particular, $L(0)$ is the right inverse of $\phi'(0)$. For $L$ we
require (as in \cite{L-Z}) the following growth condition to hold:
\begin{equation}
|L(x)y|_n\leq C\left(|x|_{\lambda n+d}|y|_d+|y|_{\lambda n
+d}\right),\label{5}
\end{equation}
for any $(x,y)\in U\times F$ and $n\geq0$, and given $d\geq0$ and
$\lambda\geq1$. 
We also do need a bound for the Taylor rest
$$
R(x,v)=\phi(x+v)-\phi(x)-\phi'(x)v.
$$
Our hipothesis will be that:
\begin{equation}
|R(x,v)|_n\leq C\left(|x|_n|v|_l^2+|v|_l|v|_n\right).\label{6}
\end{equation}
All along the proof we will use an approximation to the inverse of
$\phi'$ we will need to bound. Namely, what we need is that there is
an $m$ so that: 
\begin{equation}
|\sum_{l=0}^{\infty}(S_\theta(1-\phi'(x)L(x)))^l y|_n\leq
C(n)\theta^{m} (|x|_n |y|_d+ |y|_n).\label{7}
\end{equation}

\subsection{the lemmas}
Let $1\leq\lambda<\tau<2$, take $\tau=\frac{\lambda+2}2$. Let us
consider the sequence  $(\theta_p)_{p\in {\mathbb N}}$ defined by 
$\theta_p=2^{\tau^p}$. Observe that $\tp^\tau=\tpu$.

We want to find the solution $x$ of the equation $\phi(x)=y$, for
$y\in E$ small enough. To do so, we define the sequence
$(x_p)_p=\left(x_p(y)\right)_p$ by:
$$
\begin{array}{rcllll}
x_0=&\!\!\!\! 0& && \\
x_{p+1}=&\!\!\!\! x_p+&\!\!\!\!\Delta x_p&&&\\
           &&\!\!\!\!\Delta x_p&\!\! = S_{\tp}L(x_p)\sum_{l=0}^\infty 
           (S_{\tp}(1-\phi'(x_p)L(x_p))      )^l&\!\!\!\!z_p&\\
                                             &&&&\!\!\!\!z_p&\!=y-\phi(x_p).
\label{8}
\end{array}
$$

\begin{lemma}
Let $y\in E$ satisfy $|y|_d\leq1$. Let us assume that $|x_j|_d<1$ for
$j=1\div p$ (in order to have a well defined sequence).  Then for any 
$n\geq d$,
$$
|x_p|_n\leq K(n)\,\tp^{L(n)}|y|_n
$$
where
$$
L(n)=\frac n\lambda\frac{\lambda-1}{\tau-1}+\frac1\lambda
\frac{d+\lambda}{\tau-1}+\frac{m}{\tau-1}.$$
\end{lemma}

\demo{Proof}
From \ref{1}, we get that:
$$
|z_j|_n\leq|y|_n+|\phi(x_j)|_n\leq C(|y|_n+|x_j|_n)
$$
Thus, if $|y|_d\leq1$ and $|x_j|_d<1$, then $|z_j|_d<2C$.

On the other hand, using \ref{7}:
\begin{eqnarray*}
|\sum_{l=0}^\infty \left(S_{\tp}(1-\phi'(x_p)L(x_p))\right)^l z_p|_n
&\leq& C {\theta_p}^m(|x_p|_n |z_p|_d+|z_p|_n)\leq \\
&\leq& C {\theta_p}^m(|x_p|_n +|y|_n)
\end{eqnarray*}

Let $k=n-\frac{n-d}{\lambda}$. Then:
\begin{eqnarray*}
|\Delta x_j|_n&=&|S_{\tau_j}L(x_j)\sum_{l=0}^\infty 
           \left(S_{\tau_j}(1-\phi'(x_j)L(x_j))\right)^l z_j|_n\leq \\
&\leq& C(n,k)\theta_j^k|L(x_j)\sum_{l=0}^\infty 
           \left(S_{\tau_j}(1-\phi'(x_j)L(x_j))\right)^l z_j|_{n-k}\leq\\
&\leq& C\theta_j^k(|x_j|_{\lambda (n-k)+d}\,
|\sum_{l=0}^\infty \left(S_{\tau_j}(1-\phi'(x_j)L(x_j))\right)^l z_j|_d+ \\
&&+ |\sum_{l=0}^\infty 
\left(S_{\tau_j}(1-\phi'(x_j)L(x_j))\right)^l z_j|_{\lambda(n-k)+d})\leq\\
&\leq& C\theta_j^k (C\theta_j^m |x_j|_n+C\theta_j^m(|x_j|_n+|y|_n) )\leq\\
&\leq& C\theta_j^{k+m}(|x_j|_n+|y|_n).
\end{eqnarray*}

By using repeatedly this bound, and that $x_0=0$, we obtain:
\begin{eqnarray*}
|x_{p+1}|_n&\leq& |x_p|_n+|\Delta x_p|_n\leq\\
&\leq&( C\tp^{k+m}+1)|x_p|_n+C\tp^{k+m}|y|_n\leq\\
&\leq& (C\tp^{k+m}+1)[(C\theta_{p-1}^{k+m}+1)|x_{p-1}|_n
+C\theta_{p-1}^{k+m}|y|_n]+C\tp^{k+m}|y|_n\leq\\
&\leq&\left(\prod_{j=0}^{p}(C\theta_j^{k+m}+1)\right)|x_0|_n+
|y|_n\sum_{j=p}^{0}C\theta_j^{k+m}\prod_{i=p}^{j+1}(C\theta_i^{k+m}+1)\leq \\
&\leq& |y|_n\sum_{j=p}^{0}C\theta_j^{k+m}\prod_{i=p}^{j+1}(2C\theta_i^{k+m})\leq\\
&\leq& |y|_n\sum_{j=p}^{0}(2C)^{p+1}\prod_{i=p}^{0}2^{(k+m)\tau^i}=\\
&=&|y|_n(p+1)(2C)^{p+1}2^{(k+m)\frac{\tau^{p+1}-1}{\tau-1}}.
\end{eqnarray*}
Let us notice that
$(\tau-1)L(n)=n+1-\frac{n-d}{\lambda}+m=k+m+1>k+m$, hence:
$$
\frac{(p+1)(2C)^{p+1}\left(2^{\frac{k+m}{\tau-1}}\right)^{\tau^{p+1}-1}}
{\left(2^{L(n)}\right)^{\tau^{p+1}}}
\stackrel{p\to\infty}{\longrightarrow}0.
$$
Therefore there is a $K=K(n)$ so that
$$
|x_{p+1}|_n\leq(p+1)(2C)^{p+1}2^{(k+m)\frac{\tau^{p+1}-1}{\tau-1}}
\leq
K2^{L(n)\tau^{p+1}}=K\tpu^{L(n)}.
$$
\qed

Our next goal will be to see that 
$|z_p|_d\stackrel{p\to\infty}{\longrightarrow}0$, whenever $y\in F$ is
small enough. Namely, we will prove:

\begin{lemma}
There are $M$, $s_0$, and $\delta>0$ so that if $|y|_{s_0}<\delta$ and
$|x_j|_d<1$ for $j=0\div p$ we have the bound:
$$
|z_p|_d\leq M\tp^{-\mu}|y|_{s_0}\qquad\mbox{ where }\mu=
\frac{2+\tau}{2-\tau}(d+m).
$$
\end{lemma}

\demo{Proof}
By induction on $p$. Observe that $x_{p+1}=x_p+\Delta x_p$, so
$\phi(x_{p+1})=\phi(x_p)+\phi'(x_p)\Delta x_p+R(x_p,\Delta x_p)$, 
and hence, if we write
$A=1-\phi'(x_p)L(x_p)$, we have that:
\begin{eqnarray*}
z_{p+1}&
=&y-\phi(x_{p+1})=y-\phi(x_p)-\phi'(x_p)\Delta x_p-R(x_p,\Delta x_p)=\\
&=&z_p-\phi'(x_p)S_{\tp}L(x_p)
\sum_{l=0}^\infty \left(S_{\tp}(1-\phi'(x_p)L(x_p))\right)^l z_p
-R(x_p,\Delta x_p)=\\
&=&z_p-\phi'(x_p)L(x_p)\sum_{l=0}^\infty (S_{\tp}A)^l z_p+\\
&&+\phi'(x_p)(1-S_{\tp})L(x_p)\sum_{l=0}^\infty (S_{\tp}A)^l z_p
-R(x_p,\Delta x_p)=\\
&=&z_p+(1-\phi'(x_p)L(x_p))\sum_{l=0}^\infty (S_{\tp}A)^l z_p
-\sum_{l=0}^\infty (S_{\tp}A)^l z_p+\\
&&+\phi'(x_p)(1-S_{\tp})L(x_p)\sum_{l=0}^\infty (S_{\tp}A)^l z_p
-R(x_p,\Delta x_p)=\\
&=&z_p+\sum_{l=0}^\infty S_{\tp}A(S_{\tp}A)^l z_p+
(1-S_{\tp})A\sum_{l=0}^\infty (S_{\tp}A)^l z_p-\\
&&-\sum_{l=0}^\infty (S_{\tp}A)^l z_p
+\phi'(x_p)(1-S_{\tp})L(x_p)\sum_{l=0}^\infty (S_{\tp}A)^l z_p-\\
&&-R(x_p,\Delta x_p)=\\
&=&z_p+\sum_{l=0}^\infty (S_{\tp}A)^{l+1} z_p+
(1-S_{\tp})A\sum_{l=0}^\infty (S_{\tp}A)^l z_p
-\sum_{l=0}^\infty (S_{\tp}A)^l z_p+\\
&&+\phi'(x_p)(1-S_{\tp})L(x_p)\sum_{l=0}^\infty (S_{\tp}A)^l z_p
-R(x_p,\Delta x_p)=\\
&=&\overbrace{\phi'(x_p)(1-S_{\tp})L(x_p)\sum_{l=0}^\infty (S_{\tp}A)^l z_p}^{(a)}+
\overbrace{(1-S_{\tp})A\sum_{l=0}^\infty (S_{\tp}A)^l z_p}^{(b)}-\\
&&-\overbrace{R(x_p,\Delta x_p)}^{(d)}
\end{eqnarray*}

To bound (a), we use \ref{2}, \ref{i}, \ref{4}, and \ref{7}. 
Then, for any $s\geq d$, and $s_0=\lambda s+d$:
\begin{eqnarray*}
&&|\phi'(x_p)(1-S_{\tp})L(x_p)\sum_{l=0}^\infty (S_{\tp}A)^l z_p|_d\leq\\
&&\leq C(|x_p|_d |(1-S_{\tp})L(x_p)\sum_{l=0}^\infty (S_{\tp}A)^l z_p|_d+\\
&&+|(1-S_{\tp})L(x_p)\sum_{l=0}^\infty (S_{\tp}A)^l z_p|_d)\leq\\
&&\leq 2C |(1-S_{\tp})L(x_p)\sum_{l=0}^\infty (S_{\tp}A)^l z_p|_d \leq\\
&&\leq C \tp^{-(s-d)} |L(x_p)\sum_{l=0}^\infty (S_{\tp}A)^l z_p|_s \leq\\
&&\leq C \tp^{-(s-d)}(|x_p|_{s_0} |\sum_{l=0}^\infty (S_{\tp}A)^l z_p|_d +
|\sum_{l=0}^\infty (S_{\tp}A)^l z_p|_{s_0})\leq \\
&&\leq C \tp^{-(s-d)}(C \tp^m |x_p|_{s_0} + 
C\tp^m (|x_p|_{s_0}|z_p|_d+|z_p|_{s_0})) \leq \\
&&\leq C \tp^{-(s-d-m)}(|x_p|_{s_0} + |z_p|_{s_0}) \leq \\
&&\leq C \tp^{-(s-d-m)}(|x_p|_{s_0} + |y|_{s_0}) \leq \\
&&\leq C \tp^{-(s-d-m)}(K\tp^{L(s_0)} + 1) |y|_{s_0} \leq \\
&&\leq C \tp^{-(s-d-m-L(s_0))} |y|_{s_0}.
\end{eqnarray*}

But
\begin{eqnarray*}
s-d-m-L(s_0)&=&s-d-m-(\lambda s+d)
\frac1\lambda \frac{\lambda-1}{\tau-1}-
\frac1\lambda \frac{d+\lambda}{\tau-1}-\frac{m}{\tau-1}=\\
&=&s\frac{\tau-\lambda}{\tau-1}-d-
\frac1\lambda \frac{d+\lambda}{\tau-1}- m\frac{\tau}{\tau-1}.
\end{eqnarray*}
Thus this expression tends to infinity as $s\to \infty$, so that we
can choose $s$ so as to get $s-d-m-L(s_0)\geq \mu\tau$. 
For $\tau=\frac{2+\lambda}{2}$ we need
$s\geq\frac{P_3(d,\lambda)}{2(2-\lambda)^2}$.

For such an $s$ we have seen that:
\begin{eqnarray*}
|\phi'(x_p)(1-S_{\tp})L(x_p)\sum_{l=0}^\infty (S_{\tp}A)^l z_p|_d&\leq&
C \tp^{-(s-d-m-L(s_0))}|y|_{s_0}\leq \\
&\leq& C\tpu^{-\mu}|y|_{s_0}.
\end{eqnarray*}

To bound (b) we proceed in a similar way than as for (a), and we get that:
\begin{eqnarray*}
&&|(1-S_{\tp})A \sum_{l=0}^\infty (S_{\tp}A)^l z_p|_d \leq\\
&&\leq C \tp^{-(s_0-d)} |A \sum_{l=0}^\infty (S_{\tp}A)^l z_p|_{s_0}  \leq\\
&&\leq C \tp^{-(s_0-d)}(|x_p|_{s_0} |\sum_{l=0}^\infty (S_{\tp}A)^l z_p|_d +
|\sum_{l=0}^\infty (S_{\tp}A)^l z_p|_{s_0})\leq \\
&&\leq C \tp^{-(s_0-d)}(C \tp^m |x_p|_{s_0} + 
C\tp^m (|x_p|_{s_0}|z_p|_d+|z_p|_{s_0})) \leq \\
&&\leq C \tp^{-(s_0-d-m)}(|x_p|_{s_0} + |z_p|_{s_0}) \leq \\
&&\leq C \tp^{-(s_0-d-m)}(|x_p|_{s_0} + |y|_{s_0}) \leq \\
&&\leq C \tp^{-(s_0-d-m)}(K\tp^{L(s_0)} + 1) |y|_{s_0} \leq \\
&&\leq C \tp^{-(s_0-d-m-L(s_0))} |y|_{s_0}.
\end{eqnarray*}
Since $s_0=\lambda s+d$, this expression tends to infinity faster that
the former one. Thus for the same $s$ as before we have that:
$$
|(1-S_{\tp})A \sum_{l=0}^\infty (S_{\tp}A)^l z_p|_d\leq
C \tp^{-(s_0-d-m-L(s_0))}|y|_{s_0}\leq
C\tpu^{-\mu}|y|_{s_0}.
$$

To bound the Taylor series we use \ref{5}, getting:
\begin{eqnarray}
|R(x_p,\Delta x_p)|_d&\leq& C(|x_p|_d|\Delta x_p|_d^2+|\Delta x_p|_d^2)
\leq 2C |\Delta x_p|^2_d\leq\nonumber \\
&\leq& C (\tp^d)^2|L(x_p)\sum_{l=0}^\infty (S_{\tp}A)^l z_p|_0^2\leq\nonumber \\
&\leq& C \left(\tp^d\right)^2(|x_p|_d|\sum_{l=0}^\infty (S_{\tp}A)^l z_p|_d+
|\sum_{l=0}^\infty (S_{\tp}A)^l z_p|_d)^2\leq\nonumber \\
&\leq& C \left(\tp^d\right)^2|\sum_{l=0}^\infty (S_{\tp}A)^l z_p|_d^2\leq\nonumber \\
&\leq& C \left(\tp^{d+m}\right)^2(|x_p|_d |z_p|_d + |z_p|_d)^2\leq\nonumber \\
&\leq&  C \tp^{2d+2m}|z_p|_d^2.
\label{9}
\end{eqnarray}
Then, as $-2<-\tau$,
$$
2d+2m-2\mu= 2(d+m)\frac{-2\tau}{2-\tau}=\frac{-4}{2-\tau}(d+m)\tau<
\frac{-2-\tau}{2-\tau}(d+m)\tau =-\mu \tau
$$

Now from the induction hypothesis we get that:
$$
|R(x_p,\Delta x_p)|_d\leq C\tp^{2d+2m}|z_p|_d^2\leq
CM^2\tp^{2d+2m}\tp^{-2\mu}|y|_{s_0}^2\leq C M^2 \tpu^{-\mu}|y|_{s_0}^2.
$$

Hence what we have seen is that
$$
|z_{p+1}|_d\leq C(1+M^2|y|_{s_0})\tpu^{-\mu}|y|_{s_0},
$$
For some $C$ not deppending on $p$. WLOG, $M>C$, and we choose
$\delta$ so that  $\delta\leq \min\{1,\frac{M-C}{CM^2}\}$. 
Then the lemma is satisfied at least for 
$|y|_{s_0}<\delta$.
\qed

\begin{lemma}
There is a $\delta>0$ so that, for $|y|_{s_0}<\delta$, $|x_j|_d<1$ for
all $j\geq 0$.
\end{lemma}

\demo{Proof}
Again by induction on $j$. For $j=0$ the requirements are trivially
satisfyed. Assume that the requierements are satisfyed for a certain
$j$. Then if  $|x_j|_d<1,\, j=0\div p$, we have that $|\Delta
x_j|_d\leq C\theta_j^{d+m}|z_p|_d$, if we follow the same procedure as
in \ref{9}. Because of lemma 2, we have that $|\Delta x_j|_d\leq
C\theta_j^{-(\mu-d-m)}|y|_{s_0}$, therefore
 $|x_{p+1}|_d\leq\sum_{j=0}^p|\Delta x_j|_d\leq
C\left(\sum_{j=0}^p\theta_j^{-(\mu-d-m)}\right)|y|_{s_0}$. Using that
$\mu-d-m=(d+m)\frac{2\tau}{2-\tau}>0$, if we choose 
$\delta<\min\{\delta,\left(
\sum_{j=0}^{\infty}\theta_j^{-(\mu-d-m)}\right)^{-1}\}$, we get
\begin{equation}
|x_{p+1}|_d\leq\left(\sum_{j=0}^{\infty}\theta_j^{-(\mu-d-m)}\right)|y|_{s_0} \leq
\left(\sum_{j=0}^{\infty}\theta_j^{-(\mu-d-m)}\right)\delta<1.\label{10}
\end{equation}

Hence lemmas 1 and 2 are true for any $|y|_{s_0}\leq \delta$, with no
restrictions on the sequence $(x_p)_p$.\qed

Next we will try to improve the estimate of lemma 2, changing $\mu$ by
any $a>0$.

\begin{lemma}
For any $a>0$ there are constants $C=C(a)$ and $n=n(a)$ so that
$$
|z_p|_d\leq C |y|_{n(a)}\tp^{-a}
$$
for any  $p\geq0$ and any $y$ with $|y|_{s_0}<\delta$.
\end{lemma}

\demo{Proof}
Because of lemma 2, the statement of the lemma is satisfied for
$0<a\leq\mu$.  Let $a\geq\mu$ and assume the statement is true for
this $a$. We will see that the statement is also satisfied for  $a+d+m$.

We use that
\begin{eqnarray*}
|z_{p+1}|_d&\leq& 
\overbrace{|\phi'(x_p)(1-S_{\tp})L(x_p)
      \sum_{l=0}^\infty (S_{\tp}A)^l z_p|_d}^{(a)}+\\
&&+\overbrace{|(1-S_{\tp})A\sum_{l=0}^\infty (S_{\tp}A)^l z_p|_d}^{(b)}+
\overbrace{|R(x_p,\Delta x_p)|_d}^{(d)}
\end{eqnarray*}
and we choose $n_0=\lambda n+d$ with $n-d-m-L(n_0)>\tau(a+d+m)$. Then,
$$
|(1-S_{\tp})L(x_p) \sum_{l=0}^\infty (S_{\tp}A)^l z_p|_d 
\leq C {\theta_{p+1}}^{-(a+d+m)}|y|_{n_0}.
$$

$$
|(1-S_{\tp})A \sum_{l=0}^\infty (S_{\tp}A)^l z_p|_d 
\leq C {\theta_{p+1}}^{-(a+d+m)}|y|_{n_0}.
$$

We can bound the Taylor rest using 
$$|R(x_p,\Delta x_p)|_d\leq
C|\Delta x_p|_d^2\leq C\tp^{2d+2m}|z_p|_d^2$$
 and applying lemma 4 for this $a$ to this formula. by doing so, we get
$$
|R(x_p,\Delta x_p)|_d\leq C \tp^{-2(a-d-m)}|y|_{n(a)}^2\leq
C \tp^{-2(a-d-m)}|y|_{2n(a)}
$$
where, for the last inequality, we have used that $|y|_0\leq 1$ and
the convexity inequalities. On the other hand, 
$2a-2(d+m)\geq\tau(a+d+m)\iff a(2-\tau)\geq (d+m)(2+\tau)\iff a\geq
\mu$. 
Therefore
$|R(x_p,\Delta x_p)|_d\leq C \tpu^{-(a-d-m)}|y|_{2n(a)}$. 
If we choose 
$C=\max\{C,\theta_0^{(a-d-m)}\}$ we have  $|z_0|_d\leq |y|_d\leq
C\theta_0^{-(a-d-m)}|y|_{n(a+d+m)}$, and, using the induction, also that
$$
|z_{p+1}|_d\leq C|y|_{n(a+d+m)}\tpu^{-(a+d+m)}
$$
for any $p\geq0$, where $n(a+d+m)=\max\{2n(a),n_0\}$.
\qed

With these bounds on $|z_p|_d$ we can obtain new bounds on $|z_p|_n$,
using the convexity inequalities. To do so, we proceed in the
following way: 

\begin{lemma}
For any $n\geq0$ and any $b>0$ there are $C=C(n,b)$ and $\sigma(n,b)$
so that for any $y\in E$ with $|y|_{s_0}<\delta$ and any $p\geq0$ we
have that:
\begin{eqnarray*}
|\Delta x_p|_n&\leq&C|y|_{\sigma(n,b)}\tp^{-b}\\
|z_p|_n&\leq&C|y|_{\sigma(n,b)}\tp^{-b}
\end{eqnarray*}
\end{lemma}

\demo{Proof}
From the convexity inequalities,
$$
|\Delta x_p|_n\leq C|\Delta x_p|_0^{1/2}|\Delta x_p|_{2n}^{1/2}
$$
and, using lemma 1,
$$
|\Delta x_p|_{2n}\leq |x_{p+1}|_{2n}+
|x_p|_{2n}\leq C\left(\tpu^{L(2n)}+\tp^{L(n)}\right)|y|_{2n}\leq C
\tp^{\tau L(2n)}|y|_{2n}.
$$
From lemma 4 we know that 
$|\Delta x_p|_0\leq C{\theta_{p}}^{m}|z_p|_d\leq C|y|_{n(a)}\tp^{-a+m}$, 
for any $a>0$. By taking
$a=2b+m+\tau L(2n)$ and $\sigma(b,n)=\max\{n(a),2n\}$ what we get from
the previous inequality is that:
\begin{eqnarray*}
|\Delta x_p|_n&\leq& C
\tp^{\frac{-a+m}{2}}\tp^{\tau \frac{L(2n)}{2}}|y|_{\sigma(b,n)}=\\
&=&C\tp^{\frac{-2b-m-\tau L(2n)+m}{2}+\tau\frac{L(2n)}{2}}|y|_{\sigma(b,n)}
=C\tp^{-b}|y|_{\sigma(b,n)}
\end{eqnarray*}

Likewise, for $z_p$ we use that 
$|z_p|_n\leq C|z_p|_0^{\frac12}|z_p|_{2n}^{\frac12}$. Then we use the
definition of $z_p$ and lemma 1 obtaining that 
$|z_p|_{2n}\leq |y|_{2n}+C|x_p|_{2n} \leq C
|y|_{2n}+\tp^{L(n)}|y|_{2n}\leq C\tp^{L(n)}|y|_{2n}$. 
Hence, since the norms are increasing, we see that:
\begin{eqnarray*}
|z_p|_n&\leq& C|z_p|_d^\unmig|z_p|_{2n}^\unmig\leq
C\tp^{-\frac{a}{2}}\tp^{\frac{L(n)}{2}}|y|_{\sigma(b,n)}=\\
&=&C\tp^{-b-\frac{m}2-\frac{\tau-1}{2}L(n)}|y|_{\sigma(b,n)}
\leq C\tp^{-b}|y|_{\sigma(b,n)}
\end{eqnarray*}
with the same  $\sigma(b,n)$ as before.\qed

\begin{theorem}
Under the hypothesis  \ref{1}, \ref{2}, \ref{3}, \ref{4},
\ref{5}, \ref{6}, and \ref{7}, $\phi$ has a local right inverse $\psi$,
defined on a neighborhood $V=\{y\in F \,|\,|y|_{s_0}<\delta\}$ of $0$, and
satisfying $|\psi(y)|_d\leq C|y|_{s_0}$. 
\end{theorem}

\demo{Proof}
Using that $|x_p-x_{p+l}|_n\leq C\sum_{j=p}^{p+l-1}|\Delta x_p|_n
|y|_{\sigma(b,n)}$ and lemma 5, we get that:
$$
|x_p-x_{p+l}|_n\leq C(n,b) |y|_{\sigma(b,n)}\sum_{j=p}^\infty \theta_{j}^{-b}
$$
and thus, for $y\in V$, $(x_p)_p$ is a Cauchy sequence. Let
$x=\lim_{p\to\infty}x_p$. As, again because of lemma 5,
$z_p=y-\phi(x_p)\to 0$, and since  $\phi$ is continuous,
 $y=\lim_{p\to\infty}\phi(x_p)=\phi(x)$. For $y\in V$, we define
 $\psi(y)=\lim_{p\to\infty}x_p(y)$, and thus 
$\phi\circ\psi(y)=I$. \ref{10} gives us the desired bound.
\qed

\end{document}